\documentclass[graybox]{svmult}

\usepackage{amsmath}
\usepackage{amssymb}
\usepackage{tabu}
\usepackage{multirow}
\usepackage{url}

\usepackage{type1cm}        
\usepackage{makeidx}         
\usepackage{graphicx}        
\usepackage{multicol}        
\usepackage[bottom]{footmisc}

\usepackage{newtxtext}       %
\usepackage{newtxmath}       

\makeindex             

\newcommand{\pd}[2]{\frac{\partial#1}{\partial#2}}
\renewcommand{\div}{\nabla\cdot}
\newcommand{\grad}{\nabla\,}
\newcommand{\bi}[1]{^{(#1)}}
\newcommand{\op}[1]{\mathcal{#1}}
\newcommand{\opinv}[1]{\mathcal{#1}^{-1}}

\begin{document}

\title*{Multipreconditioning with application to two-phase incompressible Navier--Stokes flow}
\titlerunning{Multipreconditioning with application to two-phase flow}
\author{Niall Bootland and Andrew Wathen}
\institute{Niall Bootland \at University of Strathclyde, Department of Mathematics and Statistics, Glasgow, UK.\\ \email{niall.bootland@strath.ac.uk}
\and Andrew Wathen \at University of Oxford, Mathematical Institute, Oxford, UK\\ \email{andy.wathen@maths.ox.ac.uk}\\\\
This publication is based on work supported by the EPSRC Centre for Doctoral Training in Industrially Focused Mathematical Modelling (EP/L015803/1) in collaboration with the US Army Coastal and Hydraulics Laboratory and HR Wallingford. In particular, we would  like to thank Dr Chris Kees for his expert guidance in this project.}
%
%
\maketitle

\abstract*{We consider the use of multipreconditioning to solve linear systems when more than one preconditioner is available but the optimal choice is not known. In particular, we consider a selective multipreconditioned GMRES algorithm where we incorporate a weighting that allows us to prefer one preconditioner over another. Our target application lies in the simulation of incompressible two-phase flow. Since it is not always known if a preconditioner will perform well within all regimes found in a simulation, we also consider robustness of the multipreconditioning to a poorly performing preconditioner. Overall, we obtain promising results with the approach.}

\abstract{We consider the use of multipreconditioning to solve linear systems when more than one preconditioner is available but the optimal choice is not known. In particular, we consider a selective multipreconditioned GMRES algorithm where we incorporate a weighting that allows us to prefer one preconditioner over another. Our target application lies in the simulation of incompressible two-phase flow. Since it is not always known if a preconditioner will perform well within all regimes found in a simulation, we also consider robustness of the multipreconditioning to a poorly performing preconditioner. Overall, we obtain promising results with the approach.}

\section{Introduction}
\label{sec:1-Introduction}

In challenging fluid flow simulations used to model hydraulic processes it is often not clear what the best choice of preconditioner might be for solving a given linear system $\op{A}\vec{x} = \vec{b}$. Further, disparate flow regimes can be encountered in a simulation and the optimal preconditioner may change throughout. One can imagine trying to adaptively change the preconditioner based on tracking the current flow regime. However, this requires knowing \emph{a priori} which preconditioner is likely best in any given regime as well as a suitable evaluation of the current flow, which may well vary within the domain. The required sophistication and good prior knowledge of the preconditioners' performance makes such an adaptive approach less appealing.

Instead we consider using multiple preconditioners simultaneously, aiming to get the best of each. If we can combine the preconditioners then we would like to know whether we can achieve performance similar to the (unknown) best preconditioner and, further, if together they provide an improvement over any individual approach. Another key question to ask would be that of robustness: whether inclusion of a poorly performing preconditioner significantly affects the overall performance.

These ideas are encompassed within multipreconditioning strategies, where either the iterative method or preconditioning incorporates more than one preconditioner. There are several ways in which multipreconditioning can be employed but it is salient to consider the computational cost incurred weighed against the performance improvements that might be gained. Note, however, that such a strategy might not simply be aiming to give the optimal performance for solving a given system but to provide an overall robustness during a simulation spanning differing regimes.

A simple way to incorporate multiple preconditioners into an iterative method is to change the preconditioner at each iteration, in which case a flexible solver such as FGMRES \cite{Saad-FGMRES} is required. This is exemplified in \emph{cycling}, where the preconditioner choice changes in a prescribed cyclic order \cite{RuiYongAndChen}. However, results show convergence never better than the best choice of preconditioner on its own; though such a choice is unknown in advance. While only observed empirically, it stands to reason that this is unlikely to provide improvement over the best preconditioner for any given linear system, though it may help provide robustness over a sequence of problems.

Another strategy is to form a single preconditioner from the options available. This is employed in \emph{combination} preconditioning, in which the action of the inverse of the preconditioner is a linear combination of other preconditioner inverses. The term was introduced in \cite{StollAndWathen} and pursued further in \cite{PestanaAndWathen}, however, their main focus is on maintaining symmetry or positive definiteness (in some nonstandard inner product) so more efficient iterative methods can be used. Nonetheless, combination preconditioning could equally be applied to nonsymmetric cases with less restriction on requiring certain parameter choices or need for a nonstandard inner product.

A similar idea, using linear combinations of preconditioned operators, is found in the earlier \emph{multi-splitting} method \cite{O'LearyAndWhite}. The idea is to utilise multiple different splitting methods to solve the linear system. The approach can be thought of as a stationary iteration with each splitting providing a preconditioner. Yet, as with combination preconditioning, fixed 	weights for the contributions must be chosen in advance.

Except for cycling, these approaches allow for parallelism in the application of multiple preconditioners. However, the performance of the underlying iterative method will depend on the overall effectiveness of the preconditioners and how they are combined. Instead, we consider a multipreconditioned GMRES method \cite{GreifReesAndSzyld} that retains the parallelisable application of preconditioners but computes weights as part of the algorithm which are, in some sense, optimal. It considers not just one new search direction at each iteration but several, given by each preconditioner. We note that the idea was first applied to the conjugate gradient method for symmetric positive definite systems in \cite{BridsonAndChen}. However, with multiple preconditioners the search space grows exponentially fast as we continue to iterate. Thus, a selective variant of the algorithm which restricts this growth to be linear is typically necessary.

\section{Multipreconditioned GMRES (MPGMRES)}
\label{sec:2-MultipreconditionedGMRES}

In the standard preconditioned GMRES (or FGMRES) method, at each iteration a new search direction, based on the preconditioned operator, is added to the search space and then a least-squares problem is solved to find a solution with minimum residual norm. The key idea behind \emph{multipreconditioned GMRES} (MPGMRES) \cite{GreifReesAndSzyld} is to add \emph{multiple} new search directions at each iteration coming from the different preconditioners available. In fact, the method adds all new search directions from combinations of the preconditioned operators applied to vectors in the current search space, making the search space very rich. An Arnoldi-type block procedure is then used to obtain an orthonormal basis of the search space. MPGMRES then computes the optimal new iterate from this space in the minimum residual least-squares sense. Hence, note that the weights defining the contributions from each preconditioned operator are computed as part of the procedure, unlike in other approaches.

To understand how this \emph{complete} MPGMRES algorithm works, suppose we have $\ell \ge 2$ preconditioners $\op{P}_{i}$, $i = 1, \ldots, \ell$. We start with an initial residual vector $\vec{r}\bi{0}$, which we normalise to give the first basis vector $V\bi{1} = \beta^{-1} \vec{r}\bi{0}$, with $\beta = \|\vec{r}\bi{0}\|_{2}$, and collect together the preconditioned (normalised) residuals
\begin{align}
\label{MPGMRES:Z1}
Z\bi{1} = \beta^{-1} \left[\opinv{P}_{1}\vec{r}\bi{0},\ldots,\opinv{P}_{\ell}\vec{r}\bi{0}\right] \in \mathbb{R}^{n \times \ell}.
\end{align}
Using an Arnoldi-type block procedure we orthogonalise columns of $W = \op{A}Z\bi{1}$ with respect to our current basis $V\bi{1}$ and amongst themselves by using a reduced QR factorisation. Normalising then provides new basis vectors $V^{(2)} \in \mathbb{R}^{n \times \ell}$.

At each iteration, $k$, we increase the MPGMRES search space by applying each of the preconditioners to our newest basis vectors $V\bi{k}$, computing
\begin{align}
\label{MPGMRES:Zk}
Z\bi{k} = \left[\opinv{P}_{1}V\bi{k},\ldots,\opinv{P}_{\ell}V\bi{k}\right] \in \mathbb{R}^{n \times \ell^k}.
\end{align}
The Arnoldi-type block procedure is then used to orthogonalise $W = \op{A}Z\bi{k}$ with respect to the current basis $\widetilde{V}_{k} = \left[V\bi{1} \ \ldots \ V\bi{k}\right]$ and within itself. This yields new basis vectors $V\bi{k+1} \in \mathbb{R}^{n \times \ell^k}$ and, by storing the coefficients from the Arnoldi-type step in an upper Hessenberg matrix $\widetilde{H}_{k}$, we obtain an Arnoldi-type decomposition
\begin{align}
\label{MPGMRES:ArnoldiDecomposition}
\op{A}\widetilde{Z}_{k} = \widetilde{V}_{k+1}\widetilde{H}_{k},
\end{align}
where $\widetilde{Z}_{k} = \left[Z\bi{1} \ \ldots \ Z\bi{k}\right]$. Note that any linear dependency in columns of $\widetilde{Z}_{k}$, due to redundancy in the user-provided preconditioners, can be avoided using deflation; see \cite[\S3]{GreifReesAndSzyld}. Now that we have a search space then, similarly to FGMRES, we solve a linear least-squares problem for the minimum residual solution to
\begin{align}
\label{MPGMRES:MinResProb}
\min_{\vec{x} \in \vec{x}\bi{0} + \operatorname{range}\left(\widetilde{Z}_{k}\right)} \left\| \vec{b} - \op{A}\vec{x} \right\|_{2} = \min_{\vec{y}} \left\| \|\vec{r}\bi{0}\|_{2} \;\! \vec{e}_{1} - \widetilde{H}_{k}\vec{y}, \right\|_{2},
\end{align}
where $\vec{x} = \vec{x}\bi{0} + \widetilde{Z}_{k}\vec{y}$. Note that there is a natural generalisation of the standard GMRES polynomial minimisation property, as detailed in \cite{GreifReesAndSzyld}.

While the search space for complete MPGMRES is very rich, we note that it grows exponentially at each iteration, and thus becomes prohibitive in practice. As such, a variant which selects only some of the potential search directions, ideally ensuring only linear growth, is natural to consider as a more practical alternative.

\section{Selective MPGMRES (sMPGMRES)}
\label{sec:3-SelectiveMPGMRES}

To balance the benefits gained by adding multiple search directions with the storage and compute costs, we might wish to fix the number of preconditioner applications and matrix--vector products independent of the iteration, allowing for parallelisation of these operations via use of a fixed number of processors. To do so, we consider limiting the growth of the search space to be linear with respect to the iteration number $k$ by using a \emph{selective MPGMRES} (sMPGMRES) algorithm outlined in \cite{GreifReesAndSzyld}.

The search directions in MPGMRES are given by a collection of column vectors $Z$. To limit the growth of the search space we limit the size of $Z$, in particular to be proportional to the number of preconditioners, independent of $k$. To do this we select only certain search directions from the span of the columns of $Z$, giving a selective MPGMRES algorithm. There are many strategies to choose these directions, for instance, instead of applying the preconditioners to all columns of $V\bi{k}$, as in \eqref{MPGMRES:Zk}, we might apply them to just a single vector from $V\bi{k}$, selecting this vector differently for each preconditioner. This selection choice need not be the same at each iteration and could incorporate randomness if desired. The corresponding $Z\bi{k}$ is then
\begin{align}
\label{sMPGMRES:Zk-select}
Z\bi{k} = \left[\opinv{P}_{1}V\bi{k}_{:,s_{1}},\ldots,\opinv{P}_{\ell}V\bi{k}_{:,s_{\ell}}\right],
\end{align}
where $V\bi{k}_{:,s_{i}}$ is the $s_{i}$th column of $V\bi{k}$ and $s_{i}$ might change with $k$.

An alternative to applying each preconditioner to just one vector from $V\bi{k}$ is to apply them all to a linear combination of these vectors, namely to $V\bi{k}\vec{\alpha}\bi{k}$ for some vector $\vec{\alpha}\bi{k}$ of appropriate size detailing the contribution from each column of $V\bi{k}$. The corresponding $Z\bi{k}$ is then
\begin{align}
\label{sMPGMRES:Zk-lincomb}
Z\bi{k} = \left[\opinv{P}_{1}V\bi{k}\vec{\alpha}\bi{k},\ldots,\opinv{P}_{\ell}V\bi{k}\vec{\alpha}\bi{k}\right] \in \mathbb{R}^{n \times \ell}.
\end{align}
Note that a natural choice for $\vec{\alpha}\bi{k}$ is the vector $\vec{1}$, of all ones. All of these selection methods result in choosing a lower dimensional subspace of the full space and then minimising over this subspace. With these selection strategies, where we limit $Z\bi{k}$ to $\ell$ new directions each iteration, $\widetilde{V}_{k+1}$ has $k\ell+1$ basis vectors while the number of columns of $\widetilde{Z}_{k}$ is $k\ell$. Hence, the storage is proportional to $k$, as in FGMRES, as opposed to exponential in $k$, like complete MPGMRES.

Now suppose we have reason to favour one preconditioner over another and, for simplicity, that there are just two candidate preconditioners $\op{P}_{1}$ and $\op{P}_{2}$. We would like our selective approach to incorporate knowledge of which preconditioner to favour. As such, we might choose an $\vec{\alpha}\bi{k} = \vec{\alpha}$ to weight more the contributions coming from one of the preconditioners. Consider the initial steps in sMPGMRES: we start with new search directions $Z\bi{1}$ and orthogonalise them to be $V\bi{2}$
\begin{align}
\label{sMPGMRES:Z1}
Z\bi{1} = \beta^{-1} \left[\opinv{P}_{1}\vec{r}\bi{0},\opinv{P}_{2}\vec{r}\bi{0}\right] \stackrel{\mathrm{orthog.}}{\longrightarrow} V^{(2)},
\end{align}
then add search directions $Z\bi{2}$ which are orthogonalised to be $V\bi{3}$
\begin{align}
\label{sMPGMRES:Z2}
Z^{(2)} = \left[\mathcal{P}_{1}^{-1}V^{(2)}\vec{\alpha},\mathcal{P}_{2}^{-1}V^{(2)}\vec{\alpha}\right] \stackrel{\mathrm{orthog.}}{\longrightarrow} V^{(3)}.
\end{align}
So $\vec{\alpha} = (\alpha_{1},\alpha_{2})^{T}$ weighs the contributions from each of the two preconditioners as $V\bi{2}\vec{\alpha} = \alpha_{1} V\bi{2}_{:,1} + \alpha_{2} V\bi{2}_{:,2}$ and the two columns of $V\bi{2}$ come from the two different preconditioned residuals. If we let $\vec{\alpha} = (\alpha,1-\alpha)^{T}$, for some $\alpha \in (0,1)$, then the parameter $\alpha$ states how much we favour the first preconditioner, with $\alpha = \frac{1}{2}$ giving equal weighting and being equivalent to using the vectors of all ones ($\vec{\alpha}=\vec{1}$), as suggested above. Similar strategies could be used to weight contributions from more than two preconditioners.

In this weighted version of sMPGMRES the ordering of the preconditioners $\op{P}_{1},\ldots,\op{P}_{\ell}$ is important as we weight them differently. However, even with equal weighting (that is, $\vec{\alpha} = \vec{1}$) ordering is important. This more nuanced asymmetry within sMPGMRES is an aspect not mentioned in \cite{GreifReesAndSzyld}. The asymmetry comes about from the need to orthogonalise the new search directions in $Z\bi{k}$ within themselves. The contribution from the first preconditioner is allowed to be in any new direction but this direction is taken out of the contribution from subsequent preconditioners, and so on as we orthogonalise \emph{in order} the contributions from all preconditioners. This means that if the direction from the last preconditioner is mostly within the span of the preceding directions it may well contribute very little of value, despite coming from a good preconditioner when applied by itself. As a general rule then, we might value less these final search directions as the useful components may have already been taken out. This suggests taking a weighting $\vec{\alpha}$ which decreases in the components, instead of being equal, might be preferred. Nonetheless, in practice with a small number of good preconditioners, $\vec{\alpha} = \vec{1}$ might suffice to be as good. We will see that when we favour a preconditioner the ordering will matter, even if we are weighting the preconditioners in the same way. Further, ordering can still have a significant impact even when just two preconditioners are used and they are weighted equally, especially when one of the preconditioners is poorer.

\section{Numerical results for sMPGMRES}
\label{sec:4-NumericalResults}

Here we apply sMPGMRES within a two-phase incompressible Navier--Stokes flow problem. That is, to solve linear systems associated with discretisation of\begin{subequations}
\label{NSModel}
\begin{align}
\label{NSMomentum}
\rho \, \pd{\vec{u}}{t} + \rho \, \vec{u} \cdot \grad \vec{u} - \div \left( \mu \left( \grad \vec{u} + (\grad \vec{u})^{T} \right) \right) + \grad p & = \rho \, \vec{f}, \\
\label{NSIncompressibility}
\div \vec{u} & = 0,
\end{align}
\end{subequations}
for velocity $\vec{u}$ and pressure $p$ where density $\rho$ and dynamic viscosity $\mu$ are piecewise constant, representing the two phases. An important dimensionless quantity that appears is the dominating Reynolds number $Re$ over the two phases, a parameter which quantifies the ratio of inertial to viscous forces within a fluid. Our results will also exhibit how performance depends on $Re$. An auxiliary equation to describe how $\rho$ and $\mu$ vary in time with the flow is required, such as a level set equation; for the full model see \cite{BootlandBentleyKeesAndWathen}. We consider seeking the $\boldsymbol{Q}_{2}$--$\,\boldsymbol{Q}_{1}$ finite element solution using Newton iteration to treat the nonlinearity. We utilise block preconditioners, in particular those introduced in \cite{BootlandBentleyKeesAndWathen}. These are two-phase versions of the \emph{pressure convection--diffusion} (PCD) and \emph{least-squares commutator} (LSC) approaches \cite{ElmanSilvesterAndWathen}. To answer questions of robustness we further use a SIMPLE-type preconditioner, also discussed in \cite{BootlandBentleyKeesAndWathen}. We restrict our results to focus on the two preconditioner case ($\ell=2$) using \eqref{sMPGMRES:Zk-lincomb} with $\vec{\alpha}\bi{k} = (\alpha,1-\alpha)^{T}$ for some $\alpha \in (0,1)$. We follow exactly the simplified problem of a lid-driven cavity used in \cite{BootlandBentleyKeesAndWathen} along with the same implementations, as such we omit the details for brevity. The only difference is we now use sMPGMRES to solve the Newton systems via the MATLAB implementation\footnote{\scriptsize\url{www.mathworks.com/matlabcentral/fileexchange/34562-multi-preconditioned-gmres}} which accompanies \cite{GreifReesAndSzyld}.

We focus on iteration counts, as opposed to timings, since our implementation runs in serial and so does not take advantage of the inherent parallelism. Note that, when we tabulate our results using sMPGMRES, the iteration counts given in bold emphasise the best choice of weighting parameter $\alpha$ which provides the minimum number of iterations for a given pair of preconditioners. The preconditioner given on the left of a set of results is used as the first preconditioner in sMPGMRES.

\begin{table}[t]
	\centering
	\caption{Average preconditioned sMPGMRES iterations upon Newton linearisation using weighted combinations of PCD and LSC with density ratio $1.2 \times 10^{-3}$, viscosity ratio $1.8 \times 10^{-2}$ (values for air-water flow), $h = 1/64$, and varying Reynolds number $Re$ and time-step $\Delta t$.}
	\label{Table:CavityNewtonPCD&LSC_Re&Deltat}
	\tabulinesep=0.15mm
	\begin{tabu}{cc|c|ccccc|c|ccccc|c}
	     & & & \multicolumn{5}{c|}{$\alpha$ in PCD--LSC} & & \multicolumn{5}{c|}{$\alpha$ in LSC--PCD} & \\
	     $\Delta t$ & $Re$ & PCD & \ 0.9 \ & \ 0.7 \ & \ 0.5 \ & \ 0.3 \ & \ 0.1 \ & LSC & \ 0.9 \ & \ 0.7 \ & \ 0.5 \ & \ 0.3 \ & \ 0.1 \ & PCD \\
	     \hline
	     \multirow{5}{*}{$10^{-1}$} & $10$ & 16 & \textbf{14} & \textbf{14} & \textbf{14} & 16 & 18 & 18 & 15 & \textbf{14} & 15 & 15 & 21 & 16 \\
	     & $10^{1.5}$ & 16 & 14 & \textbf{13} & \textbf{13} & 14 & 16 & 17 & 14 & \textbf{13} & 14 & 15 & 18 & 16 \\
	     & $100$ & 15 & 14 & \textbf{13} & 14 & 16 & 17 & 26 & 15 & 14 & \textbf{13} & 14 & 16 & 15 \\
	     & $10^{2.5}$ & 16 & 14 & \textbf{12} & \textbf{12} & \textbf{12} & 14 & 14 & 12 & \textbf{11} & 12 & 14 & 24 & 16 \\
	     & $1000$ & 19 & 15 & \textbf{13} & 14 & 15 & 15 & 19 & \textbf{13} & \textbf{13} & \textbf{13} & 32 & 38 & 19 \\
	     \hline
	     \multirow{5}{*}{$1$} & $10$ & 19 & 18 & \textbf{17} & 18 & 21 & 23 & 24 & 19 & \textbf{18} & \textbf{18} & \textbf{18} & 21 & 19 \\
	     & $10^{1.5}$ & 21 & 19 & \textbf{18} & \textbf{18} & 20 & 22 & 23 & 19 & \textbf{18} & \textbf{18} & 19 & 22 & 21 \\
	     & $100$ & 25 & 21 & \textbf{18} & 20 & 21 & 23 & 32 & 21 & 19 & \textbf{18} & 21 & 26 & 25 \\
	     & $10^{2.5}$ & 27 & 24 & 19 & \textbf{18} & \textbf{18} & 20 & 22 & 18 & \textbf{17} & 19 & 25 & 37 & 27 \\
	     & $1000$ & 31 & 27 & \textbf{24} & \textbf{24} & \textbf{24} & 26 & 34 & 25 & \textbf{23} & 26 & 37 & 63 & 31 \\
	     \hline
	     \multirow{5}{*}{$10$} & $10$ & 20 & 19 & \textbf{18} & 19 & 22 & 23 & 25 & 20 & \textbf{19} & \textbf{19} & \textbf{19} & 22 & 20 \\
	     & $10^{1.5}$ & 24 & 21 & \textbf{20} & 21 & 23 & 25 & 26 & 22 & \textbf{20} & 21 & 22 & 27 & 24 \\
	     & $100$ & 30 & 26 & \textbf{23} & 26 & 27 & 28 & 42 & 27 & 25 & \textbf{23} & 26 & 30 & 30 \\
	     & $10^{2.5}$ & 35 & 31 & 29 & \textbf{27} & 30 & 35 & 38 & 32 & 29 & \textbf{28} & 32 & 41 & 35 \\
	     & $1000$ & 44 & 47 & \textbf{40} & 44 & 41 & 49 & 58 & 43 & \textbf{38} & 39 & 46 & 85 & 44 \\
     \end{tabu}
\end{table}

Table \ref{Table:CavityNewtonPCD&LSC_Re&Deltat} displays results for combining two-phase PCD and LSC. We see that the best iteration counts are seen towards the centre of the table, that is with a weighting parameter $\alpha$ closer to $\tfrac{1}{2}$, though we see some bias towards larger $\alpha$ for both orderings as the asymmetry of ordering might suggest. In this example most choices of $\alpha$ will provide some improvement over either of PCD or LSC individually while the best choice can allow convergence using up to 32\% fewer iterations. Note that the choice $\alpha = \tfrac{1}{2}$ typically gives iterations counts close to optimum. Given that it is not clear that we necessarily should do any better than the best preconditioner by itself, these results are quite promising and show that sMPGMRES can improve performance in terms of the number of iterations required. Furthermore, we see in this case that the performance is not particularly sensitive to $\alpha$. To examine robustness, we now include the SIMPLE-type preconditioner, a method which performs poorly here.

\begin{table}[t]
	\centering
	\caption{Average preconditioned sMPGMRES iterations upon Newton linearisation using weighted combinations of LSC and SIMPLE with density ratio $1.2 \times 10^{-3}$, viscosity ratio $1.8 \times 10^{-2}$ (values for air-water flow), $h = 1/64$, and varying Reynolds number $Re$ and time-step $\Delta t$.}
	\label{Table:CavityNewtonLSC&SIMPLE_Re&Deltat}
	\tabulinesep=0.15mm
	\begin{tabu}{cc|c|ccccc|c|ccccc|c}
	     & & & \multicolumn{5}{c|}{$\alpha$ in LSC--SIMPLE} & & \multicolumn{5}{c|}{$\alpha$ in SIMPLE--LSC} & \\
	     $\Delta t$ & $Re$ & LSC & \ 0.9 \ & \ 0.7 \ & \ 0.5 \ & \ 0.3 \ & \ 0.1 \ & SIMPLE & \ 0.9 \ & \ 0.7 \ & \ 0.5 \ & \ 0.3 \ & \ 0.1 \ & LSC \\
	     \hline
	     \multirow{5}{*}{$10^{-1}$} & $10$ & \textbf{18} & \textbf{18} & 19 & 24 & 48 & 97 & 164 & 48 & 24 & 19 & \textbf{18} & 19 & \textbf{18} \\
	     & $10^{1.5}$ & 17 & \textbf{16} & 17 & 21 & 32 & 85 & 154 & 41 & 21 & 17 & \textbf{16} & 18 & 17 \\
	     & $100$ & 26 & 23 & \textbf{22} & 27 & 35 & 51 & 131 & 38 & 26 & 23 & \textbf{22} & 24 & 26 \\
	     & $10^{2.5}$ & \textbf{14} & \textbf{14} & 15 & 18 & 27 & 35 & 116 & 33 & 19 & 15 & 15 & 15 & \textbf{14} \\
	     & $1000$ & \textbf{19} & \textbf{19} & 20 & 22 & 24 & 39 & 109 & 31 & 25 & 21 & 21 & 21 & \textbf{19} \\
	     \hline
	     \multirow{5}{*}{$1$} & $10$ & 24 & \textbf{22} & 24 & 32 & 60 & 93 & 177 & 40 & 27 & 23 & \textbf{22} & 24 & 24 \\
	     & $10^{1.5}$ & 23 & \textbf{21} & 24 & 31 & 60 & 95 & 185 & 40 & 27 & 23 & \textbf{22} & 25 & 23 \\
	     & $100$ & 32 & \textbf{30} & 31 & 40 & 70 & 103 & 188 & 60 & 36 & 31 & \textbf{30} & 31 & 32 \\
	     & $10^{2.5}$ & 22 & \textbf{20} & 21 & 26 & 50 & 109 & 190 & 46 & 25 & 21 & \textbf{20} & 22 & 22 \\
	     & $1000$ & 34 & \textbf{31} & 33 & 43 & 51 & 78 & 190 & 62 & 38 & 35 & \textbf{32} & \textbf{32} & 34 \\
	     \hline
	     \multirow{5}{*}{$10$} & $10$ & 25 & \textbf{22} & 25 & 33 & 60 & 96 & 179 & 40 & 28 & \textbf{23} & \textbf{23} & 24 & 25 \\
	     & $10^{1.5}$ & 26 & \textbf{24} & 26 & 34 & 62 & 98 & 192 & 41 & 29 & \textbf{24} & 25 & 27 & 26 \\
	     & $100$ & 42 & \textbf{38} & 39 & 47 & 69 & 104 & 207 & 64 & 42 & 37 & \textbf{36} & 39 & 42 \\
	     & $10^{2.5}$ & 38 & \textbf{33} & 35 & 42 & 82 & 125 & 233 & 54 & 38 & \textbf{33} & 34 & 37 & 38 \\
	     & $1000$ & 58 & \textbf{49} & 51 & 64 & 96 & 150 & 294 & 85 & 57 & 51 & \textbf{49} & 52 & 58 \\
     \end{tabu}
\end{table}

Table \ref{Table:CavityNewtonLSC&SIMPLE_Re&Deltat} combines the LSC and SIMPLE-type preconditioners. We see that, when LSC is used as the first preconditioner, primarily there is relatively little gained from including the SIMPLE-type approach with the best choice either being to simply use LSC or else a large $\alpha$ favouring LSC, though the best reduction in iteration counts does reach to 15\%. However, if we change the ordering to have the SIMPLE-type approach first, the picture looks slightly different. While the best iteration counts are very similar, this time any $\alpha \le \tfrac{1}{2}$ gives results comparable to LSC. This suggests that, while we do not gain much in the way of improved performance, the algorithm is still fairly robust to varying $\alpha$ so long as we do not favour the poorly performing preconditioner too strongly. This example also provides a case where, with equal weighting ($\alpha = \tfrac{1}{2}$), the ordering of the preconditioners can substantially matter, with one choice giving iteration counts that are similar or better than LSC and the other giving results that are somewhat worse than LSC. Furthermore, it is by putting the worst preconditioner first (which by the asymmetry is subtly favoured) that we obtain the better results. While at first this may sound counter-intuitive, we can make sense of this observation by considering what the selection in sMPGMRES is doing. If the good preconditioner is used first then we take this contribution away from that of the second preconditioner, likely making it even worse, then by equally weighting these we are allowing a large component of this much worse contribution to prevail. On the other hand, if the worse preconditioner is first, we remove this component from the contribution of the better preconditioner, which is unlikely to make this contribution worse and may possibly make it even better. Thus we see this latter combination is more favourable than the former, though we may not expect it to provide significantly better results than the best preconditioner by itself. We note that, in results not shown, a somewhat similar scenario occurs when combining PCD and the SIMPLE-type approach; see also \cite{Bootland} for further numerical results.

Our study show promise that sMPGMRES can combine multiple preconditioners to reduce overall iteration counts and, additionally, provide robustness in situations when one preconditioner is performing poorly. Further, weights can be incorporated to favour preconditioners and results are not particularly sensitive to any sensible choice of weights, though ordering can be important. It remains to confirm how much speed-up can be gained from sMPGMRES but initial results in \cite{GreifReesAndSzyld} are positive.

\bibliographystyle{spmpsci}
\bibliography{references}

\begin{thebibliography}{10}
\providecommand{\url}[1]{{#1}}
\providecommand{\urlprefix}{URL }
\expandafter\ifx\csname urlstyle\endcsname\relax
  \providecommand{\doi}[1]{DOI~\discretionary{}{}{}#1}\else
  \providecommand{\doi}{DOI~\discretionary{}{}{}\begingroup
  \urlstyle{rm}\Url}\fi

\bibitem{Bootland}
Bootland, N.: Scalable two-phase flow solvers.
\newblock D.Phil. thesis, University of Oxford  (2018)

\bibitem{BootlandBentleyKeesAndWathen}
Bootland, N., Bentley, A., Kees, C., Wathen, A.: Preconditioners for two-phase
  incompressible {Navier--Stokes} flow.
\newblock SIAM J. Sci. Comput. \textbf{41}(4), B843--B869 (2019)

\bibitem{BridsonAndChen}
Bridson, R., Greif, C.: A multipreconditioned conjugate gradient algorithm.
\newblock SIAM J. Matrix Anal. Appl. \textbf{27}(4), 1056--1068 (2006)

\bibitem{ElmanSilvesterAndWathen}
Elman, H.C., Silvester, D.J., Wathen, A.J.: Finite Elements and Fast Iterative
  Solvers: with Applications in Incompressible Fluid Dynamics, second edn.
\newblock Oxford University Press (2014)

\bibitem{GreifReesAndSzyld}
Greif, C., Rees, T., Szyld, D.B.: {GMRES} with multiple preconditioners.
\newblock SeMA \textbf{74}(2), 213--231 (2017)

\bibitem{O'LearyAndWhite}
O'Leary, D.P., White, R.E.: Multi-splittings of matrices and parallel solution
  of linear systems.
\newblock SIAM J. Algebraic Discret. Methods \textbf{6}(4), 630--640 (1985)

\bibitem{PestanaAndWathen}
Pestana, J., Wathen, A.J.: Combination preconditioning of saddle point systems
  for positive definiteness.
\newblock Numer. Linear Algebra Appl. \textbf{20}(5), 785--808 (2013)

\bibitem{RuiYongAndChen}
Rui, P.L., Yong, H., Chen, R.S.: Multipreconditioned {GMRES} method for
  electromagnetic wave scattering problems.
\newblock Microw. Opt. Technol. Lett. \textbf{50}(1), 150--152 (2007)

\bibitem{Saad-FGMRES}
Saad, Y.: A flexible inner-outer preconditioned {GMRES} algorithm.
\newblock SIAM J. Sci. Comput. \textbf{14}(2), 461--469 (1993)

\bibitem{StollAndWathen}
Stoll, M., Wathen, A.: Combination preconditioning and the
  {Bramble--Pasciak}$^+$ preconditioner.
\newblock SIAM J. Matrix Anal. Appl. \textbf{30}(2), 582--608 (2008)

\end{thebibliography}

\end{document}